\documentclass[11pt]{amsart}
\usepackage{amsmath,amssymb,amsthm}
\usepackage[latin1]{inputenc}
\usepackage{version,tabularx,multicol}
\usepackage{graphicx,float,psfrag}
\usepackage{stmaryrd}

\theoremstyle{plain}
\newtheorem{theo}{Theorem}[section]

\newtheorem{prop}[theo]{Proposition}
 
\newtheorem{lem}[theo]{Lemma}

\theoremstyle{remark}

 \def\beqlb{\begin{eqnarray}}\def\eeqlb{\end{eqnarray}}
 \def\beqnn{\begin{eqnarray*}}\def\eeqnn{\end{eqnarray*}}

 \def\ar{\!\!&}

 \def\qed{\hfill$\Box$\medskip}

\newcommand{\bcen}{\begin{center}}
\newcommand{\ecen}{\end{center}}
\newcommand{\bgeqn}{\begin{equation}}
\newcommand{\edeqn}{\end{equation}}

\def\l{\left}
\def\r{\right}

 \def\ar{\!\!\!&}

\begin{document}

\title[maximal out-degree of GW trees]{A note on the maximal out-degree of Galton-Watson trees}
\date{\today}

\author{Xin He}

\address{Xin He, School of Mathematical Sciences, Beijing Normal University, Beijing 100875, P.R.CHINA}

\email{hexin@bnu.edu.cn}

\begin{abstract}
In this note we consider both the local maximal out-degree and the global maximal out-degree of Galton-Watson trees.
In particular, we show that the tail of any local maximal out-degree and that of the offspring distribution are asymptotically of the same order.
However for the global maximal out-degree, this is only true in the subcritical case.
\end{abstract}



\maketitle

\section{Introduction}\label{s:i}

In a Galton-Watson tree (GW tree), we call the number of offsprings of a vertex the $\emph{out-degree}$ of this vertex.
Then naturally by $\emph{maximal out-degree}$ of the GW tree we mean the maximal number of offsprings of all the vertices in the tree.
In \cite{H14} we have proposed a new way to condition random trees, that is,
condition random trees to have large maximal out-degree. Then under this new conditioning,
we studied local convergence of critical or subcritical GW trees.
Theorem 3.3 in \cite{H14} states that in the subcritical case, tail of the offspring distribution and tail of the maximal out-degree are asymptotically equivalent, apart from the ratio $1-\mu$, where $\mu$ is the expectation of the offspring distribution.
This result seems to be interesting on its own, although in \cite{H14} it is essential for the proof of local convergence in the subcritical case.

Motivated by \cite{H14}, later in \cite{HL14} we have systematically studied the maximal jump of continuous-state branching processes (CB processes).
Note that here the maximal jump of CB processes actually corresponds to the maximal out-degree of GW trees, since
vertices of L\'{e}vy trees do not branch at the same time,
after binary branches are excluded from the consideration.
In particular, we have distinguished the \emph{local maximal jump} and the \emph{global maximal jump}, which are the maximal jump
over a finite time interval and the whole time interval, respectively.
Several exact expressions and explicit asymptotics of both the local maximal jump and the global
maximal jump of CB processes are obtained. An interesting phenomenon there is that, the tail of any local maximal jump and that of the L\'{e}vy measure are asymptotically of the same order.
However for the global maximal jump, this is only true in the subcritical case.
Naturally one wants to ask: Is there a similar phenomenon
in the setting of GW trees? The answer is yes, and in this note we present this phenomenon
in the setting of GW trees.

Now let us explain our results. Lemma \ref{generating} gives an expression of the finite-dimensional distribution of the maximal out-degrees at different generations. For the maximal out-degree at a certain generation, Proposition \ref{modone} shows that its tail and that of the offspring distribution are asymptotically of the same order. For the \emph{local maximal out-degree}, which is the maximal out-degree over a finite time interval, Theorem \ref{modlocal} shows a similar situation. Note that Theorem \ref{modlocal} can actually be proved by the method used in the proof of Theorem 3.3 in \cite{H14}, however in this note we use a slightly different method based on mean value theorem and Taylor's theorem, which we feel is more revealing, among several advantages  . Then we consider the \emph{global maximal out-degree},
which is the maximal out-degree over the whole time interval.
In Proposition \ref{sub} we prove Theorem 3.3 in \cite{H14} with our new method. The global maximal out-degree in the critical case
is more subtle than that in the subcritical case. However we can still use our new method to get some partial results in Theorem \ref{critical}.

Although our papers \cite{H14,HL14} have motivated this note directly,
it has also been influenced by two interesting papers \cite{B11,B13} from Bertoin.
In fact, \cite{B11,B13} have also influenced \cite{H14,HL14} in certain results and arguments. The topics in \cite{B11,B13} are the global maximal out-degree in two special critical cases and its scaling limit.
Our setting in this note is more general, and we consider both the local maximal jump and the global maximal jump.
However we do not consider the scaling limit here, which seems to be the main motivation of \cite{B11,B13}.

In the remaining of this introduction, we review some notations of GW trees.
In the next section, we state and prove our results.

We use $\tau(p)$ to denote a GW tree with the {\em offspring distribution} $p=(p_0,p_1,p_2,\cdots)$ on non-negative integers.
We say that $p$ is $\emph{bounded}$, if the set $\{r; p_r > 0\}$ is bounded.
Use $\mathbf{P}$ to denote the underlying probability, and $\mu$ the expectation of $p$.
We assume that $\mu\in(0,\infty)$, so that $p_0<1$.
We use \emph{(sub)critical} to mean critical or subcritical, that is, $\mu\leq 1$.
It is well-known that $\tau(p)$ has finite vertices a.s. if it is (sub)critical, and has infinite vertices with positive probability if it is supercritical. For this result and more, refer to the standard reference \cite{AN72}.
Following the terminology of Galton-Watson processes, we also say that $\tau(p)$ is \emph{extinct} if it has finite vertices, and is \emph{non-extinct} if it has infinite vertices.

At time (or generation) 0, the GW tree $\tau(p)$ has a single vertex, the so-called {\em root}.
At generation 0, the root branches into several offsprings according to the offspring distribution $p$.
We use $M_n$ to denote the maximal number of offsprings of all vertices at generation $n$, and $M_{[0,n]}$ the maximal number of offsprings of all vertices in the first $n+1$ generations, including generation 0. Then clearly $\mathbf{P}\l[M_0=r\r] =p_r$. We also use $M=M_{[0,\infty)}$ to denote the supremum of  all out-degrees of $\tau(p)$.
Clearly in the (sub)critical case $M$ is finite a.s. and it is actually the maximum of all out-degrees.
Use $F(r)$ to denote the distribution function of $p$, $\bar F(r)$ the tail function.
Similarly we use $H(r)$ to denote the distribution function of $M$, $\bar H(r)$ the tail function. Note that in the supercritical case,
it is possible for $M$ to be infinite.

We use $G(x)$ to denote the generating function of $p$ ,that is, for $x\in[0,1]$,
$$
G(x)=p_0+p_1 x+\cdots+p_n x^n+\cdots.
$$
Clearly $G$ has derivatives of all orders, and
$$
\frac{\partial G(x)}{\partial x}|_{x=1-}=\mu.
$$
We write $G^2(x)=GG(x)$ for $G(G(x))$. By induction, clearly
$$
\frac{\partial G^n(x)}{\partial x}|_{x=1-}=\mu^n.
$$
Finally for any nonnegative integer $r$, we use $G_r(x)$ to denote the generating function of $p$ truncated at $r$, that is, for $x\in[0,1]$,
$$
G_r(x)=p_0+p_1 x+\cdots+p_r x^r.
$$

\section{The results}

We begin with an exact expression of the finite-dimensional distribution of the maximal out-degrees at different generations.
This result is very elementary, however here we give a general version with our applications later in mind.

\begin{lem} \label{generating} Recall that $M_i$ is the maximal out-degree at generation $i$. Then
$$
\mathbf{P}\l[M_i\leq r_i,0\leq i \leq n\r] = G_{r_0}G_{r_1}\cdots G_{r_n}(1).
$$
\end{lem}
\proof For $n=0$, trivially
$$
\mathbf{P}\l[M_0\leq r_0\r] =p_0+p_1+\cdots+p_{r_0}=G_{r_0}(1).
$$
For $n=1$, by considering the out-degree of the root we get
$$
\mathbf{P}\l[M_0\leq r_0,M_1\leq r_1\r]=p_0+p_1 \mathbf{P}\l[M_0\leq r_1\r] +\cdots+p_{r_0}\mathbf{P}\l[M_0\leq r_1\r]^{r_0}=G_{r_0}G_{r_1}(1).
$$
The general case follows by induction and a similar identity
$$
\mathbf{P}\l[M_i\leq r_i,0\leq i \leq n\r] = p_0+p_1 \mathbf{P}\l[M_{i-1}\leq r_i,1\leq i \leq n\r]+\cdots+p_{r_0}\mathbf{P}\l[M_{i-1}\leq r_i,1\leq i \leq n\r]^{r_0}.
$$

\qed

This lemma implies that
$$
\mathbf{P}\l[M_n\leq r \r] = G^n G_r(1)=G^n (p_0+\cdots+p_r).
$$
In general, the generating function $G$ has no explicit expression.
However, we can get a simple asymptotic result of $\mathbf{P}\l[M_n=r \r]$ as $r\rightarrow\infty$.

\begin{prop} \label{modone} It is always true that
$$
\mathbf{P}\l[M_n= r \r] \leq \mu^n p_r.
$$
Assume that $p$ is unbounded, then as $r\rightarrow\infty$,
$$
\mathbf{P}\l[M_n= r \r] \sim \mu^n p_r \quad \text{and} \quad \mathbf{P}\l[M_n> r \r] \sim  \mu^n\bar{F}(r),
$$
where the first limit is understood along the infinite subsequence $\{r; p_r > 0\}$.
\end{prop}
\proof The second limit is implied by the first one. For the first limit, by mean value theorem we have
$$
\mathbf{P}_x\l[M_n = r \r]= G^n (p_0+\cdots+p_r) - G^n (p_0+\cdots+p_{r-1})
= p_r \frac{\partial G^n (y)}{\partial y}|_{y=a_r},
$$
where $a_r$ is a certain point in the open interval $(p_0+\cdots+p_{r-1},p_0+\cdots+p_r)$. The inequality is obvious since
$$\frac{\partial G^n (y)}{\partial y}|_{y=a_r}\leq \frac{\partial G^n (y)}{\partial y}|_{y=1-}=\mu^n.$$
For the limit,
clearly we have as $r\rightarrow\infty$,
$$
\mathbf{P}\l[M_n= r \r]/p_r \rightarrow \frac{\partial G^n (y)}{\partial y}|_{y=1-}=\mu^n.
$$
\qed

Then we study the local maximal jump over the time interval $[0,n]$, that is, $M_{[0,n]}$.
Proposition \ref{modone} clearly implies that
\beqlb\label{i}
\mathbf{P}\l[M_{[0,n]}= r \r] \leq \sum_{i=0}^n\mathbf{P}\l[M_i= r \r]\leq p_r\sum_{i=0}^n \mu^i.
\eeqlb
We will show that $\mathbf{P}\l[M_{[0,n]}= r \r]$ and this upper bound are asymptotically the same, as in Proposition \ref{modone}.

\begin{theo} \label{modlocal}
Assume that $p$ is unbounded, then as $r\rightarrow\infty$,
$$
\mathbf{P}\l[M_{[0,n]}= r \r] \sim (1+\mu+\cdots+\mu^n)p_r \quad \text{and} \quad \mathbf{P}\l[M_{[0,n]}> r \r] \sim (1+\mu+\cdots+\mu^n)\bar{F}(r),
$$
where the first limit is understood along the infinite subsequence $\{r; p_r > 0\}$.
\end{theo}

\proof The second limit is implied by the first one, so we only need to prove the first one.
When $n=0$, trivially$\mathbf{P}\l[M_0= r \r]=p_r$.
For $n>0$, by Lemma \ref{generating} and mean value theorem, we have
\beqnn
\mathbf{P}\l[M_{[0,n]} = r \r]\ar=\ar G_r G_r^n(1) - G_{r-1} G_{r-1}^n(1)\cr
\ar=\ar G_{r-1} G_r^n(1) - G_{r-1} G_{r-1}^n(1)+p_r [G_r^n(1)]^r\cr
\ar=\ar \mathbf{P}\l[M_{[0,n-1]}= r \r]\frac{\partial G_{r-1} (y)}{\partial y}|_{y=a_r}+p_r [G_r(1)]^r,\cr
\eeqnn
where $a_r$ is a certain point in the open interval $(G_{r-1}^n(1),G_r^n(1))$.
From (\ref{i}) we know that for any $n\geq 0$, $M_{[0,n]}$ has finite expectation.
Consequently for any $n\geq 1$, $\lim_{r\rightarrow\infty}[G_r^n(1)]^r=1$, since $G_r^n(1)=\mathbf{P}\l[M_{[0,n-1]} \leq r \r]$.
It is also easy to see that for any $R$,
$$
\mu^R=\lim_{r\rightarrow\infty}\frac{\partial G_{R}(y)}{\partial y}|_{y=a_r}\leq \lim_{r\rightarrow\infty}\frac{\partial G_{r-1}(y)}{\partial y}|_{y=a_r}
\leq\lim_{r\rightarrow\infty}\frac{\partial G(y)}{\partial y}|_{y=a_r}=\mu,
$$
where $\mu^R=\sum_{i=1}^R i p_i$. Notice that $\lim_{R\rightarrow\infty} \mu^R=\mu$.
Finally assuming the statement is true for $M_{[0,n-1]}$ gives that as $r\rightarrow\infty$,
$$
\mathbf{P}\l[M_{[0,n]} = r \r]\sim \mathbf{P}\l[M_{[0,n-1]}= r \r]\mu +p_r \sim (\mu^n+\cdots+\mu+1)p_r.
$$
\qed

Note that the above theorem can also be proved by the method used in the proof of Theorem 3.3 in \cite{H14}, which we call {\em the barehands method}.
Here we used a slightly different method, which we call {\em the calculus method}.
We feel that
the most important advantage of the calculus method is that it correctly identifies $\mu$ as the derivative of the generating function $G$
at $1-$. The calculus method is also more automatic and more efficient, we believe. Now let us use this calculus method to get a generalization of both Proposition \ref{modone} and Theorem \ref{modlocal}.

\begin{prop}
Assume that $p$ is unbounded. Consider a sequence of $m$ positive integers $(n_i,1\leq i\leq m)$ such that $n_i< n_j$ for $1\leq i<j\leq m$,
and another sequence of $m$ positive integers $(r_{n_i},1\leq i\leq m)$. Let $r=\min_{1\leq i\leq m}r_{n_i}$.
Assume additionally that, $\limsup_{r\rightarrow\infty}r_j/\min_{j< i\leq m}r_{n_i}<\infty$ for any $1\leq j <n$.
Then as $r\rightarrow\infty$,
$$
\mathbf{P}\l[\bigcup_{i=1}^m \{M_{n_i}= r_{n_i},M_{n_j}< r_{n_j};1\leq j\leq m,j\neq i\}\r] \sim \sum_{i=1}^m \mu^{n_i}p_{r_{n_i}},
$$
where the limit is understood along the infinite subset $\{(r_{n_i},1\leq i\leq m);\sum_{i=1}^m p_{r_{n_i}}> 0\}$.
\end{prop}

\proof We argue by induction. When $m=1$, this is exactly Proposition \ref{modone}. When $m=2$, we need to prove that as $r=r_{n_1}\wedge r_{n_2}\rightarrow\infty$,
$$
\mathbf{P}\l[\{M_{n_1}= r_{n_1},M_{n_2}< r_{n_2}\} \bigcup \{M_{n_2}= r_{n_2},M_{n_1}< r_{n_1}\} \r]\sim\mu^{n_1} p_{r_{n_1}}+\mu^{n_2} p_{r_{n_2}}.
$$
Lemma \ref{generating} shows that the l.h.s of the above expression is
$$
G^{n_1}G_{r_{n_1}}G^{n_2-n_1-1}G_{r_{n_2}}(1)-G^{n_1}G_{r_{n_1}-1}G^{n_2-n_1-1}G_{r_{n_2}-1}(1).
$$
Then as in Theorem \ref{modlocal}, we get that as $r\rightarrow\infty$,
\beqnn
\ar\ar G^{n_1}G_{r_{n_1}}G^{n_2-n_1-1}G_{r_{n_2}}(1)- G^{n_1}G_{r_{n_1}-1}G^{n_2-n_1-1}G_{r_{n_2}-1}(1)\cr
\ar\sim\ar \mu^{n_1} \l(G_{r_{n_1}}G^{n_2-n_1-1}G_{r_{n_2}}(1)-G_{r_{n_1}-1}G^{n_2-n_1-1}G_{r_{n_2}-1}(1)\r)\cr
\ar=\ar \mu^{n_1} \l(G_{r_{n_1}-1}G^{n_2-n_1-1}G_{r_{n_2}}(1)-G_{r_{n_1}-1}G^{n_2-n_1-1}G_{r_{n_2}-1}(1)\r)\cr
\ar\ar +\mu^{n_1} p_{r_{n_1}}\l[G^{n_2-n_1-1}G_{r_{n_2}}(1)\r]^{r_{n_1}}\cr
\ar\sim\ar \mu^{n_1+1} \l(G^{n_2-n_1-1}G_{r_{n_2}}(1)-G^{n_2-n_1-1}G_{r_{n_2}-1}(1)\r)+\mu^{n_1} p_{r_{n_1}}\cr
\ar\sim\ar \mu^{n_2} p_{r_{n_2}}+\mu^{n_1} p_{r_{n_1}},\cr
\eeqnn
where in the last step we used the induction for $m=1$.
Clearly when $m>1$, we can always reduce the case of $m$ into the case of $m-1$ in this way.

Notice that we then only need to check that $\lim_{r\rightarrow\infty}\l[G^{n_2-n_1-1}G_{r_{n_2}}(1)\r]^{r_{n_1}}=1$. Recall from the proof of Theorem \ref{modlocal} that
$\lim_{r_{n_2}\rightarrow\infty}[G_{r_{n_2}}^{n_2-n_1}(1)]^{r_{n_2}}=1$.
Finally $\lim_{r\rightarrow\infty}[G_{r_{n_2}}^{n_2-n_1}(1)]^{r_{n_1}}=1$ by our assumption, so $\lim_{r\rightarrow\infty}\l[G^{n_2-n_1-1}G_{r_{n_2}}(1)\r]^{r_{n_1}}=1$.
The general case can be proved similarly.
\qed

From now on we study the global maximal out-degree $M=M_{[0,\infty)}$.
We begin with a characterization of the distribution function of $M$.
In the (sub)critical case this result
has appeared in \cite{B13}, see page 794-795 of \cite{B13}. Here we mainly add the result in the supercritical case.

\begin{lem} \label{modglobal} Assume that $p$ is supercritical with $p_0+\cdots+p_r<1$ or (sub)critical, then
$$
\mathbf{P}\l[M \leq r \r] =q_{[0,r]},
$$
where $q_{[0,r]}$ is the unique solution of the equation $G_r(t)=t$ for $t\in[0,1]$.

\end{lem}

\proof By considering the out-degree of the root we see that
$$
\mathbf{P}\l[M \leq r \r]=p_0+p_1 \mathbf{P}\l[M \leq r \r] +\cdots+p_r\l(\mathbf{P}\l[M \leq r \r]\r)^r.
$$
So $\mathbf{P}\l[M \leq r \r]$ is a solution of the equation $G_r(t)=t$ for $t\in[0,1]$.

Then we only need to show the uniqueness of the solution. If $p_1+2p_2\cdots+rp_r\leq 1$,
then it is easy to see that $G'_r(t)< 1$ for $t\in[0,1)$. Since trivially $G_r(0)=p_0\geq 0$ and $G_r(1)=p_0+\cdots+p_r\leq 1$,
there is indeed a unique solution of the equation $G_r(t)=t$ for $t\in[0,1]$.

If $p_1+2p_2\cdots+rp_r> 1$,
then clearly $p$ is supercritical and $p_0+p_1<1$. Let $q$ be the smallest solution of $G(t)=t$ for $t\in[0,1]$ (see e.g., Section I.3 in \cite{AN72}). Trivially $G_r(t)\leq G(t)$ and $G'_r(t)\leq G'(t)$. So $G'_r(t)\leq G'(t)\leq G'(q)<1$ for $t\leq q$ (again, see Section I.3 in \cite{AN72}).
Since trivially $G_r(0)=p_0\geq 0$ and $G_r(q)\leq G(q)=q$, we see that there is a unique solution of the equation $G_r(t)=t$ for $t\in[0,q]$. Since $G_r(t)\leq G(t)$, so if $G_r(t)=t$ has a solution for $t\in(q,1]$, then $G_r(1)=1=p_0+\cdots+p_r$. We are done.
\qed

Of course, If $p$ is supercritical with $p_0+\cdots+p_r=1$, then trivially $\mathbf{P}\l[M \leq r \r] =1$.
In this case, we may just define $q_{[0,r]}=1$.
Or we may just say that $q_{[0,r]}$ is the largest solution of the equation $G_r(t)=t$ for $t\in[0,1]$.
Let $q_r=q_{[0,r]}-q_{[0,r-1]}$ for $r\geq 1$, and $q_0=q_{[0,0]}$.
Then $\mathbf{P}\l[M =r \r] =q_r$.

From now on we study the tail of $M$. We first treat the subcritical case. Note that this result was first proved in Theorem 3.3 of \cite{H14}, by the barehands method.
Here we prove it by the calculus method, which we feel is more revealing. Recall that $H(r)$ is the distribution function of $M$ and
$\bar H(r)$ the tail function, so $H(r)=q_{[0,r]}$. Also recall the simple fact that if $p_0>0$, then $q_r>0$ if and only if
$p_r>0$, see Lemma 3.1 in \cite{H14}.

\begin{prop} \label{sub} Assume that $p$ is subcritical and unbounded, then
$$
\lim_{r\rightarrow\infty}\frac{\bar{F}(r)}{\bar{H}(r)}=\lim_{r\rightarrow\infty}\frac{p_r}{q_r}=1-\mu,
$$
where the last limit is understood along the infinite subsequence $\{r; p_r > 0\}$.
\end{prop}

\proof The first limit is implied by the second one. For the second limit,
by Lemma \ref{modglobal} we have
$$
q_{[0,r]}=G_r(q_{[0,r]})=G_{r-1}(q_{[0,r]})+p_r q_{[0,r]}^r, \quad q_{[0,r-1]}=G_{r-1}(q_{[0,r-1]}).
$$
By considering the difference of the these two identities and mean value theorem, we get
\beqlb\label{subcritical}
q_r =q_r \frac{\partial G_{r-1}(y)}{\partial y}|_{y=a_r}+p_rq_{[0,r]}^r,
\eeqlb
where $a_r$ is a certain point in the open interval $(q_{[0,r-1]},q_{[0,r]})$.
Similar to (\ref{i}), we have
$$
\mathbf{P}\l[M= r \r] \leq \mathbf{P}\l[M_0= r \r]+\mathbf{P}\l[M_1= r \r]+\cdots\leq \frac{p_r}{1-\mu}.
$$
So $M$ has finite expectation, and consequently $\lim_{r\rightarrow\infty}q_{[0,r]}^r=1$.
We are done since in Theorem \ref{modlocal} we have proved that
\beqlb\label{noproblem}
\lim_{r\rightarrow\infty}\frac{\partial G_{r-1}(y)}{\partial y}|_{y=a_r}=\mu.
\eeqlb
\qed

The supercritical case is essentially trivial, since it can be reduced to the subcritical case.
Let $q$ be the probability of extinction, that is, the smallest solution of $G(t)=t$ for $t\in[0,1]$ (see e.g., Section I.3 in \cite{AN72}).
By a well-known result on supercritical GW trees (see e.g., Theorem I.12.3 on page 52 of \cite{AN72}), for any supercritical and unbounded offspring distribution $p$ with $p_0>0$, we have $0<q<1$,
$$
\lim_{r\rightarrow\infty}\bar{H}(r)=1-q, \quad \text{and}\quad \lim_{r\rightarrow\infty}\frac{p_r q^{r-1}}{q_r}=1-\sum_{i\geq 1}ip_iq^{i-1},
$$
where the last limit is understood along the infinite subsequence $\{r; p_r > 0\}$. If $p_0=0$, then trivially $q_\infty=1$.

Finally we treat the critical case.
Theorem \ref{modlocal} implies that in the critical case, $\bar{H}(r)$ is asymptotically larger than $\bar{F}(r)$. To prove this, just notice that for any $n$,
$$
\liminf_{r\rightarrow\infty}\frac{\bar{H}(r)}{\bar{F}(r)}=\liminf_{r\rightarrow\infty}\frac{\mathbf{P}\l[M> r \r] }{\bar{F}(r)}\geq \lim_{r\rightarrow\infty}\frac{\mathbf{P}\l[M_{[0,n]}> r \r] }{\bar{F}(r)}=n+1.
$$
Apart from this simple observation, generally speaking the critical case is much more subtle than the subcritical or supercritical case, as shown by Bertoin
in \cite{B11,B13}. However we can still use our calculus method to derive several partial results. We first present a simple lemma on $q_{[r,\infty)}$.

\begin{lem}\label{width} Assume that $p$ is (sub)critical, then $q_{[r,\infty)}\leq 1/r$  for any positive integer $r$.
\end{lem}

\proof This follows from a simple inequality on width of GW trees. Use $X_n$ to denote the total number of nodes of the tree $\tau$ at generation $n$ and define the width $W(\tau)$ of the tree $\tau$ by $W(\tau)=\sup_{n\geq 0} X_n$. Now by optional sampling for supermartingales, we get $\mathbf{P}[W\geq r]\leq 1/r$ for any positive integer $r$, which in turn implies the desired inequality on $q_{[r,\infty)}$.
\qed

Note that a part of the following theorem was first appeared in Lemma 1 of \cite{B13}.

\begin{theo}\label{critical} Assume that $p$ is critical and unbounded, then for $q_r$ we have
$$
\lim_{r\rightarrow\infty}\frac{q_r}{p_r}=\infty,
 \quad
q_r \leq \frac{p_r}{\sum_{i>r}(ip_i)},
 \quad\text{and} \quad
\limsup_{r\rightarrow\infty}\frac{q_r^2}{p_r}\leq \frac{2}{\sigma^2},
$$
where $\sigma^2$ is the variance of $p$ (could be infinite) and the limits are understood along the infinite subsequence $\{r; p_r > 0\}$. For $\bar{H}(r)$ similarly we have
$$
\lim_{r\rightarrow\infty}\frac{\bar{H}(r)}{\bar{F}(r)}=\infty,
 \quad
\bar{H}(r) \leq \frac{\bar F(r)}{\sum_{i>r}(ip_i)},
 \quad\text{and} \quad
\limsup_{r\rightarrow\infty}\frac{\bar{H}^2(r)}{\bar{F}(r)}\leq \frac{2}{\sigma^2}.
$$
\end{theo}

\proof For the first statement, recall (\ref{subcritical}) and (\ref{noproblem}). By Lemma \ref{width}, we see that $\liminf_{r\rightarrow\infty} q_{[0,r]}^r\geq 1/e$. We are done with the first statement since $1-\mu=0$ for a critical offspring distribution $p$.

For the second and third statements, apply Taylor's theorem to the difference of $H(r-1)=G_{r-1}(H(r-1))$ and $H(r)=G_r(H(r))=G_{r-1}(H(r))+p_rH^r(r)$ to get
$$
q_r =q_r \frac{\partial G_{r-1}(y)}{\partial y}|_{y=1-}-\frac{1}{2} q_r^2\frac{\partial^2 G_{r-1}(y)}{\partial y^2}|_{y=a_r}+p_rH^r(r),
$$
where $a_r$ is a certain point in the open interval $(H(r-1),H(r))$. So
$$
q_r\leq \frac{p_r}{1- \frac{\partial G_r(y)}{\partial y}|_{y=1-}}= \frac{p_r}{\sum_{i>r}(ip_i)}
$$
and
$$
\limsup_{r\rightarrow\infty}\frac{q_r^2}{p_r}\leq \limsup_{r\rightarrow\infty}\frac{2}{\frac{\partial^2 G_{r-1}(y)}{\partial y^2}|_{y=a_r}}=\frac{2}{\sigma^2}.
$$

For the last two statements, again apply Taylor's theorem to the difference of $H(r)=G_r(H(r))$ and $1=G_r(1)+\bar F(r)$ to get
$$
\bar H(r) =\bar H(r) \frac{\partial G_r(y)}{\partial y}|_{y=1-}-\frac{1}{2}\bar H^2(r)\frac{\partial^2 G_r(y)}{\partial y^2}|_{y=a_r}+\bar F(r),
$$
where $a_r$ is a certain point in the open interval $(H(r),1)$. So
$$
\bar{H}(r) \leq \frac{\bar F(r)}{1- \frac{\partial G_r(y)}{\partial y}|_{y=1-}}= \frac{\bar F(r)}{\sum_{i>r}(ip_i)}
$$
and
$$
\limsup_{r\rightarrow\infty}\frac{\bar H^2(r)}{\bar F(r)}\leq \limsup_{r\rightarrow\infty}\frac{2}{\frac{\partial^2 G_r(y)}{\partial y^2}|_{y=a_r}}=\frac{2}{\sigma^2}.
$$
\qed

In particular, if the critical offspring distribution $p$ has tail $\bar F(r)\sim r^{-2}$ as $r\rightarrow\infty$, then Theorem \ref{critical} shows that
$\bar H(r)= o\,(r^{-1})$ as $r\rightarrow\infty$.


\end{document}